\documentclass[a4paper,10pt]{article}
\usepackage{amsmath}
\usepackage{amsmath,amssymb,latexsym,color}
\usepackage[mathscr]{eucal}
\usepackage{mathrsfs}
\usepackage{graphpap}
\newtheorem{thm}{Theorem}[section]
\newtheorem{prop}[thm]{Proposition}
\newtheorem{lemma}[thm]{Lemma}
\newtheorem{cor}[thm]{Corollary}
\newtheorem{example}[thm]{Example}
\newtheorem{defin}[thm]{Definition}
\newtheorem{re}{Remark}
\newtheorem{co}{Comment}

\newcommand{\proof}{{\it Proof.\quad}}
\newcommand{\qed}{\hfill\Box\medskip}

\usepackage{CJK}

\begin{document}
\begin{CJK*}{GBK}{song}

\renewcommand{\baselinestretch}{1.3}
\title{More on the Terwilliger algebra of   Johnson schemes}

\author{
Benjian Lv\textsuperscript{a}\quad Carolina  Maldonado\textsuperscript{b}\quad Kaishun Wang\textsuperscript{a}\footnote{Corresponding author. E-mail address: wangks@bnu.edu.cn}\\
{\footnotesize  \textsuperscript{a} \em  Sch. Math. Sci. {\rm \&}
Lab. Math. Com. Sys.,
Beijing Normal University, Beijing, 100875,  China}\\
{\footnotesize  \textsuperscript{b} \em FCEFyN Universidad
Nacional de C\'ordoba, CIEM-CONICET, Argentina} }
 \date{}
 \maketitle

\begin{abstract}
In [F. Levstein, C. Maldonado, The Terwilliger algebra of the
Johnson schemes, Discrete Math. 307 (2007) 1621--1635], the
Terwilliger algebra of the Johnson scheme $J(n,d)$ was determined
when $n\geq 3d$. In this paper, we determine the Terwilliger algebra ${\mathcal T}$ for the remaining case
$2d\leq n<3d$.

\medskip


\end{abstract}


\section{Introduction}

Let ${\mathcal X}=(X, \{R_i\}_{i=0}^{d} )$ denote a commutative
association scheme, where $X$ is a finite set. Suppose ${\rm
Mat}_X(\mathbb C)$ denotes the   algebra over $\mathbb C$ consisting
of all matrices whose rows and columns are indexed by $X.$ For each
$i$, let $A_i$ denote the binary matrix in ${\rm Mat}_X(\mathbb C)$
whose $(x, y)$-entry is $1$ if and only if $(x, y) \in R_i.$ We call
$A_i$   the \emph{$i$th adjacency matrix} of $\mathcal X.$ We
abbreviate $A=A_1$, and call it the \emph{adjacency matrix} of
$\mathcal X.$ The subalgebra of ${\rm Mat}_X(\mathbb C)$ spanned by
$A_0,A_1,\ldots,A_d $  is called the \emph{Bose-Mesner algebra} of
$\mathcal X$, denoted by $\mathcal B.$ Since $\mathcal B$ is
commutative and generated by real symmetric matrices, it has a basis
consisting of primitive idempotents, denoted by $E_0=\frac 1{|X|}
J,E_1,E_2,\ldots,E_d$. For each $i\in \{0,1,\ldots,d\},$   write
$$
A_i=\sum_{j=0}^{d}p_i(j)E_j,\quad E_i=\frac{1}{|X|}\sum_{j=0}^{d}q_i(j)A_j.
$$
The scalars $p_i(j)$ and $q_i(j)$ are called the {\em eigenvalues}
and the {\em dual eigenvalues} of ${\mathcal X},$ respectively.

Fix $x\in X.$ For $0\leq i\leq d$, let $E^*_i$ denote the diagonal
matrix in ${\rm Mat}_X(\mathbb{C})$ whose $(y, y)$-entry is defined
by
$$
(E^*_i)_{yy} =\left\{
 \begin{array}{cl}
    1, & {\rm if}\ (x, y) \in R_i,\\
    0, & {\rm otherwise}.
   \end{array}
   \right.
$$
The subalgebra $\mathcal{T}$ of ${\rm Mat}_X(\mathbb{C})$ generated
by $A_0,A_1,\ldots, A_d; E^*_0,E^*_1,\ldots,E^*_d $ is called the
\emph{Terwilliger algebra} of $\mathcal X$ with respect to $x.$

Terwilliger \cite{Terwilliger1} first introduced the Terwilliger
algebra of association schemes, which is an important tool in
considering the structure of an association scheme.   For more
information; see
\cite{Caughman1,Caughman2,Terwilliger2,Terwilliger3}. The
Terwilliger algebra is a finite-dimensional semisimple $\mathbb
C$-algebra,   it is difficult to determine its structure in general.
The structures of the Terwilliger algebras of some association
schemes were determined; see \cite{Balmaceda, Bannai1} for group
schemes, \cite{Tomiyama} for   strongly regular graphs, \cite{Go,
Levstein1} for Hamming schemes, \cite{Levstein} for Johnson schemes,
\cite{Kong} for odd graphs, \cite{klw} for  incidence graphs of
Johnson geometry.

Let $[n]$ denote the set $\{1,2,\ldots, n\}$ and ${[n]\choose d}$
denote the collection of all $d$-element subsets of $[n]$.  For
$0\leq i\leq d$, define $R_i =\{(x, y)\in {[n]\choose d} \times
{[n]\choose d} \mid |x\cap y|=d-i\}$. Then
 $({[n]\choose
d}, \{R_i\}_{i=0}^{d} )$ is a symmetric association scheme of
 class $d$, which is called  the \emph{Johnson scheme}, denoted
by $J(n, d).$ Note that  $J(n, d)$  is isomorphic  to  $J(n, n-d)$.
So we always assume that $n\geq 2d.$

In \cite{Levstein}, the
Terwilliger algebra of the Johnson scheme $J(n,d)$ was determined when $n\geq 3d$.
In this paper, we focus on the remaining case,  and determine the
Terwilliger algebra $\mathcal T$ of   $J(n,d)$. In Section 2, we
introduce   intersection matrices and some useful identities. In
Sections 3,
  two families of subalgebras  $\mathcal{M}^{(n, d)}$  and $\mathcal{N}$
of ${\rm Mat}_X(\mathbb{C})$ are constructed. In the last two
sections, we show that $\mathcal T=\mathcal{M}^{(n, d)}$ and
$\mathcal T=\mathcal{N}$ according to $n>2d$ and $n=2d$,
respectively.

\section{Intersection matrix}

In this section we first introduce some useful identities for
intersection matrices,    then describe  the adjacency matrix of the
Johnson scheme $J(n, d)$ in terms of intersection matrices.

  Let $H^r_{i,j}(v)$ be a binary
matrix whose rows and columns are indexed by the elements of
${[v]\choose i}$ and ${[v]\choose j}$ respectively,  whose
$\alpha_i\alpha_j$-entry is $1$ if and only if $ |\alpha_i \cap
\alpha_j |= r$. For simplicity,   write  $H_{i,j}:=
H^{\min{(i,j)}}_{i,j}$. Now we introduce some useful identities for
intersection matrices.

\begin{lemma}{\rm(\cite[Theorem 3]{Mohammad-Noori})}\label{cy1}
For $0\leq l\leq \min(i,j)$ and $0\leq s\leq \min(j,k),$
$$
H^l_{i,j}(v)H^s_{j,k}(v)=\sum^{\min(i,k)}_{g=0}\Big{(}\sum^g_{h=0}{g\choose
h}{i-g\choose l-h}{k-g\choose s-h}{v+g-i-k\choose
j+h-l-s}\Big{)}H^g_{i,k}(v).
$$
\end{lemma}

\begin{lemma}{\rm(\cite[Lemma 4.4]{Levstein})}\label{hh}
Let $v$ be a positive integer.

{\rm(i)} For $0\leq i\leq j\leq l\leq v,$
$$
H_{i,j}(v)H_{j,l}(v)={l-i\choose l-j}H_{i,l}.
$$

 {\rm(ii)} For
$0\leq\max(i,l)\leq j\leq v,$
$$
H_{i,j}(v)H_{j,l}(v)=\sum_{m=0}^{j-\max(i,l)}{v-\max(i,l)-m\choose
j-\max(i,l)-m}H_{i,l}^{\min(i,l-m)}(v).
$$

{\rm(iii)} For $0\leq j\leq \min(i,l)\leq v,$
$$
H_{i,j}(v)H_{j,l}(v)=\sum_{m=0}^{\min(i,l)-j}{\min(i,l)-m\choose
j}H_{i,l}^{\min(i,l-m)}(v).
$$
\end{lemma}

Pick $x\in {[n]\choose d}$. For $0\leq i\leq d$, write $\Omega_i :=
\{y\in  {[n]\choose d}\mid |x\cap y| = d-i\}.$ Then we have the
partition ${[n]\choose d} =\dot{\cup}^d_{i=0}\Omega_i.$ Now we
consider the $m$th adjacent matrix $A_m$ of $J(n, d)$ as a block
matrix with respect to this partition. Denote
$(A_m)_{|\Omega_i\times\Omega_j}$ the submatrix of $A_m$ with rows
indexed by $\Omega_i$ and columns indexed by $\Omega_j.$

In the remaining of this paper, we always assume that $I^{(v,k) }$
denotes the identity matrix of size $v\choose k$ and $A^{ (v, k)
}_{m}$ denotes the $m$-th adjacency matrix of $J(v, k).$ In fact
$A^{ (v, k) }_{m}=H^m_{k,k}(v).$
\begin{lemma}{\rm(\cite{Levstein})}\label{str1}
 For $0\leq
i<j\leq d,$ we have
\begin{eqnarray*}
&&A_{|\Omega_i\times\Omega_i} = I^{(d,d-i)}\otimes A^{(n-d,i)} +
A^{(d,d-i)}\otimes I^{(n-d,i)},\\
&&A_{|\Omega_i\times\Omega_{i+1}}=H_{d-i,d-i-1}(d)\otimes H_{i,i+1}(n-d),\\
&&A_{|\Omega_i\times\Omega_j}=0, \ \ \ {\rm if}\ j\geq i+2,
\end{eqnarray*}
where $``\otimes"$ denotes the Kronecker product of matrices.

\end{lemma}

\section{Two algebras}

Let $2d\leq n$ and  $X={[n]\choose d}$. In this section we shall
construct two subalgebras  $\mathcal{M}^{(n, d)}$  and $\mathcal{N}$
of ${\rm Mat}_X(\mathbb{C})$, which
 are in fact  the Terwilliger algebras $\mathcal{T}$ of $J(n,d)$ according to $2d< n<
3d$ and $J(2d,d)$.

 Let $\mathcal{B}^{(v,k)}$ denote the Bose-Mesner
algebra of $J(v,k)$, and $\{E^{(v,k)}_r\}_{r=0}^{\min(k,v-k)}$
denote its primitive idempotents. Let $p_i^{(v,k)}(j)$ be the
eigenvalue of $A_i^{(v,k)}$ satisfying
$A_i^{(v,k)}E_j^{(v,k)}=p_i^{(v,k)}(j)E_j^{(v,k)}.$

For $ i, j=0,1,\ldots,   d,$  define the vector space
$$
M_{i,j}^{(n,d)}=(\mathcal{B}^{(d,d-i)}H_{d-i,d-j}(d))\otimes
(\mathcal{B}^{(n-d,i)}H_{i,j}(n-d)).
$$
Let $L: Y\longmapsto L(Y)$ be a linear map from $M_{i,j}^{(n,d)}$ to
${\rm Mat}_X(\mathbb{C})$, where $(L(Y ))_{\Omega_l\times\Omega_m}$
is $Y$ if $l = i$ and $m = j;$ $0$ otherwise. Define
\begin{eqnarray}\label{defm}
\mathcal{M}^{(n,d)}&=&\bigoplus^d_{i,j=0} L(M^{(n,d)}_{i,j}).
\end{eqnarray}

\begin{lemma}{\rm(\cite[Lemma 4.4]{Levstein})}\label{b1}
Let $R(v,k,h) := \{r \mid H^r_{k,h}(v) \neq 0\}.$ Then
$\mathcal{B}^{(v,k)}H_{k,h}(v)=H_{k,h}(v)\mathcal{B}^{(v,h)}=
\langle\{H^r_{k,h}(v)\}_{r\in R(v,k,h)}\rangle.$
\end{lemma}

 Lemma~\ref{b1} implies that
$H_{d-i,d-j}^r(d)\otimes
H_{i,j}^s(n-d),\quad {r\in R(d,d-i,d-j),s\in R(n-d,i,j)}
$ is a basis of
$M^{(n,d)}_{i,j}. $  By Lemma~\ref{cy1}, we observe that
$\mathcal{M}^{(n,d)}$ is an algebra.

By Lemmas~\ref{str1} and ~\ref{b1}, the adjacency matrix  $A$ of
$J(n, d)$ belongs to $\mathcal{M}^{(n,d)}$. Since each $m$th
adjacency matrix of $J(n, d)$ may be written as a polynomial of $A$,
one gets $\mathcal B^{(n, d)}\subseteq\mathcal{M}^{(n,d)}$.  The
fact that $E^*_m=L(H_{d-m,d-m}(d)\otimes H_{m,m}(n-d))\in
\mathcal{M}^{(n,d)}$ implies that $\mathcal{T}$ is a subalgebra of
$\mathcal{M}^{(n,d)}.$

Next we   construct another algebra $\mathcal{N}.$ For $ i,
j=0,1,\ldots,   d,$ let $N_{i,j} $ be the vector space generated by
\begin{equation}\label{nij}
(E_r^{(d,d-i)}\otimes E_s^{(d,i)}+E_s^{(d,d-i)}\otimes
E_r^{(d,i)})(H_{d-i,d-j}(d)\otimes H_{i,j}(d)),\quad 0\leq r\leq
s\leq\min(d,d-i).
\end{equation}
Define
\begin{equation}\label{n}
\mathcal{N}=\bigoplus^d_{i,j=0}L({N}_{i,j}).
\end{equation}
Observe that $\mathcal N\subseteq  \mathcal{M}^{(2d,d)}.$

\begin{lemma} Each vector space $N_{i,j}$ has the basis
\begin{equation}\label{basis}
H_{d-i,d-j}^{d-i-j+g}(d)\otimes
H_{i,j}^{h}(d)+H_{d-i,d-j}^{d-i-j+h}(d)\otimes H_{i,j}^{g}(d), \quad
g,h\in R(d,i,j),\; g\leq h.\end{equation}
\end{lemma}
\proof
By Lemma \ref{b1} we have that $\mathcal{B}^{(v,k)}H_{k,h}(v)=H_{k,h}(v)\mathcal{B}^{(v,h)}=
\langle\{H^r_{k,h}(v)\}_{r\in R(v,k,h)}\rangle.$

Also note that it holds
\begin{eqnarray*}
H_{d-i,d-j}^{d-i-j+g}(d)\otimes
H_{i,j}^{l}(d)+H_{d-i,d-j}^{d-i-j+l}(d)\otimes H_{i,j}^{g}(d)
=H_{i,j}^{g}(d)\otimes H_{i,j}^{l}(d)+H_{i,j}^{l}(d)\otimes
H_{i,j}^{g}(d),
\end{eqnarray*}
and that $\{H_{i,j}^{g}(d)\otimes
H_{i,j}^{h}(d)\}_{g,h\in R(i,j)}$ are linearly independent, which implies that\\
$\{H_{i,j}^{g}(d)\otimes H_{i,j}^{h}(d)+H_{i,j}^{h}(d)\otimes
H_{i,j}^{g}(d)\}_{g, h\in R(d,i,j),g\leq h}$ are linearly
independent. Hence, the desired result follows.
 $\qed$

\begin{thm}
Let $\mathcal{N}$ be as in {\rm (\ref{n})}. Then $\mathcal{N}$  is
an algebra.
\end{thm}
\proof It suffices to show that $N_{i,j}N_{j,k}\subseteq N_{i,k}.$
By Lemma~\ref{cy1}  there exist scalars $\alpha^g_{l,s}$ such that
$$
H_{i,j}^{l}(d)H_{j,k}^{s}(d)=\sum_{g=0}^{\min(i,k)}\alpha^g_{l,s}H_{i,k}^{g}(d),
$$
which implies that
\begin{eqnarray*} &&(H_{i,j}^{m}(d)\otimes
H_{i,j}^{n}(d)+H_{i,j}^{n}(d)\otimes H_{i,j}^{m}(d))
(H_{j,k}^{s}(d)\otimes
H_{j,k}^{t}(d)+H_{j,k}^{t}(d)\otimes H_{j,k}^{s}(d))\\
&=&\sum_{g=0}^{\min(i,k)}\sum_{l=0}^{\min(i,k)}(\alpha^g_{m,s}\alpha^l_{n,t}+
\alpha^g_{m,t}\alpha^l_{n,s})(H_{i,k}^{g}(d)\otimes
H_{i,k}^{l}(d)+H_{i,k}^{l}(d)\otimes H_{i,k}^{g}(d)),
\end{eqnarray*}
as desired. $\qed$

\begin{lemma}\label{stru2d} Let $\mathcal T$ be the Terwilliger algebra of
$J(2d,d).$ Then $\mathcal{T}$ is a subalgebra of $\mathcal{N}.$
\end{lemma}
 \proof Observe that $A_{|\Omega_i\times\Omega_{i+l}}=0$ for $l\geq2.$  By Lemma
 \ref{str1}, we obtain
 $A_{|\Omega_i\times\Omega_{i+1}}=H_{d-i,d-i-1}(d)\otimes H_{i,i+1}(d),$
 and
\begin{eqnarray*}
A_{|\Omega_i\times\Omega_i} &=& I^{(d,d-i)}\otimes
A^{(d,i)} + A^{(d,d-i)}\otimes I^{(d,i)}\\
&=&\sum\limits_{r=0}^{\min(i,d-i)}\sum\limits_{s=0}^{\min(i,d-i)}
(p_1^{(d,i)}(r)+p_1^{(d,i)}(s))(E^{(d,d-i)}_r\otimes E^{(d,i)}_s)\\
&=&\sum\limits_{r=0}^{\min(i,d-i)}\sum\limits_{s=0}^{\min(i,d-i)}\frac{p_1^{(d,i)}(r)+p_1^{(d,i)}(s)}{2}
(E^{(d,d-i)}_r\otimes E^{(d,i)}_s+E^{(d,d-i)}_s\otimes E^{(d,i)}_r).
\end{eqnarray*}
Then $A\in\mathcal{N}$. Also we have that
$E^*_m=L(H_{d-m,d-m}(d)\otimes H_{m,m}(d))\in \mathcal{N},$ so
$\mathcal{T}\subseteq\mathcal{N}.$ $\qed$

We next introduce two mappings proposed in \cite{Levstein}: the lift map denoted by $\mathcal{L}_i$ and the pullback map denoted by $\mathcal{P}_i$.

For $0\leq i<d,$ define $\mathcal{L}_i$ be the linear mapping from
$M^{(n,d)}_{i,i}$ to $M^{(n,d)}_{i+1,i+1}$ satisfying
\begin{eqnarray*}
&&\mathcal{L}_i(E_r^{(d,d-i)}\otimes E_s^{(n-d,i)})\\
&=& (H_{d-i-1,d-i}(d)E_r^{(d,d-i)}H_{d-i,d-i-1}(d))\otimes
(H_{i+1,i}(n-d)E_s^{(n-d, i)}H_{i,i+1}(n-d));
\end{eqnarray*}
for $0< i\leq d$, define $\mathcal{P}_i$ be the  linear mapping from
$M^{(n,d)}_{i,i}$ to $M^{(n,d)}_{i-1,i-1}$ satisfying
\begin{eqnarray*}
&&\mathcal{P}_i(E_r^{(d,d-i)}\otimes E_s^{(n-d,i)})\\
&=& (H_{d-i+1,d-i}(d)E_r^{(d,d-i)}H_{d-i,d-i+1}(d))\otimes
(H_{i-1,i}(n-d)E_s^{(n-d,i)}H_{i,i-1}(n-d)).
\end{eqnarray*}
Note that the lift map is defined premultiplying by $
\left(H_{d-i-1,d-i}(d)\otimes H_{i+1,i}(n-d)\right)$ (that is equal to $A_{|\Omega_{i+1}\times\Omega_{i}}$ by Lemma \ref{str1}) and post multiplying by
$\left(H_{d-i,d-i-1}(d)\otimes H_{i,i+1}(n-d)\right)$ (that is equal to $A_{|\Omega_{i}\times\Omega_{i+1}}$ by Lemma \ref{str1}). Since they belong to
$\mathcal{T}_{|\Omega_{i+1}\times\Omega_{i}}$ and $\mathcal{T}_{|\Omega_{i}\times\Omega_{i+1}}$ respectively and since $\mathcal{T}$ is an algebra, then $\mathcal{L}_i(Y)\in\mathcal{T}_{|\Omega_{i+1}\times\Omega_{i+1}}$ for any $Y\in\mathcal{T}_{|\Omega_i\times\Omega_i}.$ Similarly $\mathcal{P}_i(Y)\in\mathcal{T}_{|\Omega_{i-1}\times\Omega_{i-1}}$
for any $Y\in\mathcal{T}_{|\Omega_i\times\Omega_i}.$

\vspace{.3em}

By \cite[p.220]{Bannai}, for $0\leq j\leq k,$
\begin{eqnarray}\label{Bannai}
p_1^{(v,k)}(j)=(k-j)(v-k-j)-j ,\quad p_k^{(v,k)}(j)=(-1)^j{v-k-j\choose
k-j}.
\end{eqnarray}
Write $l_{v,k,j}=v-k+p_1^{(v,k)}(j)$ and
$p_{v,k,j}=k+p_1^{(v,k)}(j).$

\begin{lemma}{\rm(\cite[Lemma 5.6]{Levstein})}\label{llpp}
\begin{eqnarray*}
\mathcal{L}_i(E_r^{(d,d-i)}\otimes
E_s^{(n-d,i)})=p_{d,d-i,r}l_{n-d,i,s}E_r^{(d,d-i-1)}\otimes
E_s^{(n-d,i+1)},\\
\mathcal{P}_i(E_r^{(d,d-i)}\otimes
E_s^{(n-d,i)})=l_{d,d-i,r}p_{n-d,i,s}E_r^{(d, d-i+1)}\otimes
E_s^{(n-d,i-1)}.
\end{eqnarray*}
Here, $E^{(v,k)}_j=0 $ if $j>\min(k,v-k)$.
\end{lemma}

\begin{cor}{\rm \cite[Corollary 5.7]{Levstein}}\label{change}
$$
(E^{(d,d-i)}_rH_{d-i,d-j}(d))\otimes (E^{(n-d,i)}_sH_{i,j}(n-d))=
(H_{d-i,d-j}(d)E^{(d,d-j)}_r)\otimes (H_{i,j}(n-d)E^{(n-d,j)}_s).
$$
\end{cor}

\section{$\mathcal{T}$-algebra of $J(n,d)$ for $2d<n<3d$}
In this section  we always assume that    $2d<n<3d$ and $\mathcal
M^{(n, d)}$ is as in (\ref{defm}). We shall prove that
$\mathcal{M}^{(n,d)}$ is the Terwilliger algebra $\mathcal{T}$ of
$J(n,d)$.

For  any real number $a$, we have
\begin{eqnarray}\label{Aa}
(A+aI)_{|\Omega_i\times\Omega_i}=\sum^{\min(i,d-i)}_{r=0}
\sum^{\min(i,n-d-i)}_{s=0}(\mu_{i,r}+\lambda_{i,s}+a)E_r^{(d,
d-i)}\otimes E_s^{(n-d, i)},
\end{eqnarray}
where $\mu_{i,r}=p_1^{(d,d-i)}(r)$ and
$\lambda_{i,s}=p_1^{(n-d,i)}(s).$ We always assume that $a$ is a
real number large enough such that the coefficients   in (\ref{Aa})
are positive.

\begin{re}\label{V}
If these coefficients are pairwise distinct,  each
$E_r^{(d, d-i)}\otimes E_s^{(n-d, i)}$ belongs to ${\mathcal
T}_{|\Omega_i\times\Omega_i}$ since the left hand side of (\ref{Aa}) and its powers belong to ${\mathcal
T}_{|\Omega_i\times\Omega_i}$. By orthogonality of the idempotents
\begin{eqnarray}\label{Ab}
(A+aI)_{|\Omega_i\times\Omega_i}^j=\sum^{\min(i,d-i)}_{r=0}
\sum^{\min(i,n-d-i)}_{s=0}(\mu_{i,r}+\lambda_{i,s}+a)^j E_r^{(d,
d-i)}\otimes E_s^{(n-d, i)},
\end{eqnarray} obtaining a linear system of equations given by the  powers of
$(A+aI)_{|\Omega_{i} \times \Omega_{i}}$ as linear combinations of
$E_r^{(d,
d-i)}\otimes E_s^{(n-d, i)}$ with a
Vandermonde matrix. See Section 5.1 in \cite{Levstein}
for more explanation.
\end{re}

\begin{thm} Suppose $2d< n <3d$. Let $\mathcal{T}$ be the Terwilliger algebra of  $J (n,
d)$
 and $\mathcal{M}^{(n,d)}$ be the algebra as in
{\rm(\ref{defm})}. Then $\mathcal{T}=\mathcal{M}^{(n,d)}.$
\end{thm}
\proof Since
$L(M^{(n,d)}_{i,j})=L(M^{(n,d)}_{i,i}(H_{d-i,d-j}(d)\otimes
H_{i,j}(n-d))),$ by \cite[Proposition 5.2]{Levstein} it is
sufficient to prove that, for $0\leq r\leq \min(d-i,i)$ and $0\leq
s\leq \min(n-d-i,i)$,
\begin{eqnarray}\label{prove1}
E_r^{(d,d-i)}\otimes E_s^{(n-d, i)}\in \mathcal{T}_{|\Omega_i\times
\Omega_i}.
\end{eqnarray}
We shall prove (\ref{prove1}) by induction on $i$ ($i$ decreases
from $d$ to $0$).

\begin{co} Since it holds $\mathcal{P}_{i+1}(Y)\in\mathcal{T}_{|\Omega_{i}\times\Omega_{i}}$
for any $Y\in\mathcal{T}_{|\Omega_{i+1}\times\Omega_{i+1}}$ the strategy of the proof is to pull it back those projectors $E_r^{(d,d-i-1)}\otimes E_s^{(n-d, i+1)}\in \mathcal{T}_{|\Omega_i\times
\Omega_{i+1}} $ (whenever its pullback is different from zero) or separate the projectors as we explained in the Remark \ref{V}.
\end{co}

Induction: observe that
$
(A+aI)_{|\Omega_d\times\Omega_d} =
\sum_{s=0}^{d}(\mu_{d,0}+\lambda_{d,s}+a)E_0^{(d, 0)}\otimes
E_s^{(n-d,d)}.
$
Since $\lambda_{d,s}$'s are pairwise distinct, (\ref{prove1}) holds
for $i=d.$

\medskip
\textbf{Case 1.} $\lceil\frac{n-d}{2}\rceil+1\leq i\leq d-1$.

Note that in these case $\min(d-i,i)=d-i$ and $\min(n-d-i,i)=n-d-i.$

For $0\leq r\leq d-i-1$ and $0\leq s\leq n-d-i-1,$ by (\ref{Bannai})
one gets $l_{d,d-i-1,r}\neq0$ and $p_{n-d,i+1,s}\neq0$. By
Lemma~\ref{llpp},
$$
E_r^{(d,d-i)}\otimes
E_s^{(n-d,i)}=l^{-1}_{d,d-i-1,r}p^{-1}_{n-d,i+1,s}\mathcal{P}_{i+1}(E_r^{(d,d-i-1)}\otimes
E_s^{(n-d,i+1)});
$$
and so $E_r^{(d,d-i)}\otimes E_s^{(n-d,i)}\in
\mathcal{T}_{|\Omega_i\times\Omega_i}.$

It is not possible pull back
$(E_{d-i}^{(d,d-i-1)}\otimes
E_s^{(n-d,i+1)})$ not also  $(E_r^{(d,d-i-1)}\otimes
E_{n-d-i}^{(n-d,i+1)}).$

By (\ref{Aa}), we have

$\begin{array}{l}
\quad(A+aI)_{|\Omega_i\times\Omega_i}-\sum^{d-i-1}\limits_{r=0}\sum^{n-d-i-1}\limits_{s=0}(\mu_{i,r}+\lambda_{i,s}+a)E_r^{(d,
d-i)}\otimes E_s^{(n-d, i)}\\
=\sum\limits^{n-d-i}_{s=0}(\mu_{i,d-i}+\lambda_{i,s}+a)E_{d-i}^{(d,
d-i)}\otimes
E_s^{(n-d,i)}+\sum\limits^{d-i-1}_{r=0}(\mu_{i,r}+\lambda_{i,n-d-i}+a)E_{r}^{(d,d-i)}\otimes
E_{n-d-i}^{(n-d,i)}.
\end{array}$

It follows that the right hand side of the equality belongs to
$\mathcal{T}_{|\Omega_i\times\Omega_i}.$  In order to show that
(\ref{prove1}) holds, it suffices to show that each term  belongs to
$\mathcal{T}_{|\Omega_i\times\Omega_i}.$ Observe that there do not
exist three coefficients with the same value. If there exists a term
whose coefficient is different from other coefficients, then this
term belongs to $\mathcal{T}_{|\Omega_i\times\Omega_i}.$ Next
suppose that there exist two terms with the same coefficient.
Suppose that
$\mu_{i,d-i}+\lambda_{i,q}+a=\mu_{i,u}+\lambda_{i,n-d-i}+a.$ Then
$E_{d-i}^{(d, d-i)}\otimes E_q^{(n-d,i)}+E_{u}^{(d,d-i)}\otimes
E_{n-d-i}^{(n-d,i)}$ belongs to
$\mathcal{T}_{|\Omega_i\times\Omega_i},$ and by Lemma~\ref{llpp} its
image under $\mathcal P_i$ is
\begin{equation*} (i+\mu_{i,d-i})(i+\lambda_{i,q})E_{d-i}^{(d, d-i+1)}\otimes
E_q^{(n-d,i-1)}
+(i+\mu_{i,u})(i+\lambda_{i,n-d-i})E_{u}^{(d,d-i+1)}\otimes
E_{n-d-i}^{(n-d,i-1)}.
\end{equation*}

Suppose
$(i+\mu_{i,d-i})(i+\lambda_{i,q})=(i+\mu_{i,u})(i+\lambda_{i,n-d-i}).$
Since $\mu_{i,d-i}+\lambda_{i,q}=\mu_{i,u}+\lambda_{i,n-d-i},$ one
gets $\mu_{i,d-i}\lambda_{i,q}=\mu_{i,u}\lambda_{i,n-d-i}.$ It
follows that
$(\mu_{i,d-i}-\lambda_{i,n-d-i})(\mu_{i,d-i}-\mu_{i,u})=0,$ a
contradiction to $\mu_{i,d-i}\neq\lambda_{i,n-d-i} $ and
$\mu_{i,d-i}\neq\mu_{i,u}$. Therefore, we have
$(i+\mu_{i,d-i})(i+\lambda_{i,q})\neq(i+\mu_{i,u})(i+\lambda_{i,n-d-i}),$
which implies that both  $E_{d-i}^{(d, d-i+1)}\otimes
E_q^{(n-d,i-1)}$ and $E_{u}^{(d,d-i+1)}\otimes
E_{n-d-i}^{(n-d,i-1)}$ belong to $
\mathcal{T}_{|\Omega_{i-1}\times\Omega_{i-1}}.$ Computing their
image under $\mathcal L_{i-1}$, by Lemma~\ref{llpp} again
$E_{d-i}^{(d, d-i)}\otimes E_q^{(n-d,i)}$ and
$E_{u}^{(d,d-i)}\otimes E_{n-d-i}^{(n-d,i)}$ belong to $
\mathcal{T}_{|\Omega_i\times\Omega_i}$, as desired.

\medskip\textbf{Case 2.} $i=\lceil\frac{n-d}{2}\rceil.$

We divide our discussion into two subcases.

\medskip{\bf Case 2.1.} $n-d$ is odd. By Lemma~\ref{llpp}, for any
$E_r^{(d,d-i)}\otimes E_s^{(n-d,i)}\in M^{(n,d)}_{i,i},$
$$\mathcal{P}_{i}(E_r^{(d,d-i)}\otimes
E_s^{(n-d,i)})=l_{d,d-i,r}p_{n-d,i,s}E_r^{(d, d-i+1)}\otimes
E_s^{(n-d,i-1)}\neq0.
$$

Similar to the proof in Case 1, (\ref{prove1}) holds.

\medskip{\bf Case 2.2.} $n-d$ is even. For $0\leq r\leq d-i-1$ and $0\leq s\leq n-d-i-1,$ by
Lemma~\ref{llpp},
$$
E_r^{(d,d-i)}\otimes
E_s^{(n-d,i)}=l^{-1}_{d,d-i-1,r}p^{-1}_{n-d,i+1,s}\mathcal{P}_{i+1}(E_r^{(d,d-i-1)}\otimes
E_s^{(n-d,i+1)}),
$$
which implies that $E_r^{(d,d-i)}\otimes E_s^{(n-d,i)}\in
\mathcal{T}_{|\Omega_i\times\Omega_i}.$

Next we consider   $r=d-i$ or $s=n-d-i.$ Write
$\mu'_{i,r}=p_{d-i}^{(d,d-i)}(r)$ and
$\lambda'_{i,s}=p_{i}^{(n-d,i)}(s).$ Since
$(A_d)_{|\Omega_i\times\Omega_i}=A^{(d,d-i)}_{d-i}\otimes
A^{(n-d,i)}_i,$ \begin{eqnarray*}
&&(A_d+aI)_{|\Omega_i\times\Omega_i}-\sum^{d-i-1}_{r=0}\sum^{n-d-i-1}_{s=0}(\mu'_{i,r}\lambda'_{i,s}+a)E_r^{(d,
d-i)}\otimes E_s^{(n-d, i)}\\
&=&\sum^{n-d-i}_{s=0}(\mu'_{i,d-i}\lambda'_{i,s}+a)E_{d-i}^{(d,
d-i)}\otimes
E_s^{(n-d,i)}+\sum^{d-i-1}_{r=0}(\mu'_{i,r}\lambda'_{i,n-d-i}+a)E_{r}^{(d,d-i)}\otimes
E_{n-d-i}^{(n-d,i)}.
\end{eqnarray*}
It follows that the right hand side of the equality belongs to
$\mathcal{T}_{|\Omega_i\times\Omega_i}.$ By (\ref{Bannai}), observe
that $\mu'_{i,r}\lambda'_{i,n-d-i}+a$ is not equal to any other
coefficient. Then $ E_{r}^{(d,d-i)}\otimes E_{n-d-i}^{(n-d,i)}\in
\mathcal{T}_{|\Omega_i\times\Omega_i}$ for $0\leq r\leq d-i-1.$

By (\ref{Aa}),
$$
\sum^{n-d-i}_{s=0}(\mu_{i,d-i}+\lambda_{i,s}+a)E_{d-i}^{(d,
d-i)}\otimes E_s^{(n-d,i)}
$$
belongs to $\mathcal{T}_{|\Omega_i\times\Omega_i}.$ Moreover, its
coefficients are pairwise distinct,  so $E_{d-i}^{(d, d-i)}\otimes
E_s^{(n-d,i)}$ belongs to $\mathcal{T}_{|\Omega_i\times\Omega_i}$
for $0\leq s\leq n-d-i.$ Therefore, (\ref{prove1}) holds.

\ \ \newline \textbf{Case 3.} $\lceil\frac{d}{2}\rceil\leq i\leq
\lceil\frac{n-d}{2}\rceil-1.$

Note that in these case $\min(d-i,i)=d-i$ and $\min(n-d-i,i)=i.$

Similarly, by Lemma~\ref{llpp} we have that  $E_r^{(d,d-i)}\otimes
E_s^{(n-d,i)}$ belongs to $\mathcal{T}_{|\Omega_i\times\Omega_i}$
for $0\leq r\leq d-i-1$ and $0\leq s\leq i.$

By (\ref{Aa}) again,
the matrix
$$\sum^{i}_{s=0}(\mu_{i,d-i}+\lambda_{i,s}+a)E_{d-i}^{(d,
d-i)}\otimes E_s^{(n-d,i)}$$
 belongs to
$T_{|\Omega_i\times\Omega_i}.$ Moreover, its coefficients are
pairwise distinct,  so (\ref{prove1}) holds.

\ \ \newline \textbf{Case 4.} $0\leq i\leq
\lceil\frac{d}{2}\rceil-1.$

 Note that in these case $\min(d-i,i)=i$ and $\min(n-d-i,i)=i.$

By Lemma~\ref{llpp} again,  (\ref{prove1}) holds. $\qed$

Next we shall decompose $\mathcal T$ as a  direct sum of some simple
ideals.

For $0\leq r\leq \lfloor\frac{d}{2}\rfloor $ and $0\leq s\leq
\lfloor\frac{n-d}{2}\rfloor,$ define
$$
\begin{array}{rcl}
e_{r,s}&=&\min\{i\mid 0\neq E^{(d,d-i)}_r\otimes E^{(n-d,i)}_s\in
M^{(n,d)}_{i,i}\},\\
d_{r,s}&=&|\{i\mid 0\neq E^{(d,d-i)}_r\otimes E^{(n-d,i)}_s\in
M^{(n,d)}_{i,i}\}|-1.
\end{array}
$$
Note that $e_{r,s}=\max(r,s)$ and $e_{r,s}+d_{r,s}=\min(d-r,n-d-s).$

For $0\leq r\leq \min(d-i,i) $ and $0\leq s\leq \min(n-d-i,i),$
define
\begin{eqnarray}
{}^{rs}T_{ij}&=&(E^{(d,d-i)}_rH_{d-i,d-j}(d))\otimes
(E^{(n-d,i)}_sH_{i,j}(n-d)),\label{rsta}\\
{}^{rs}\mathcal{T}&=&\langle\{L(^{rs}T_{i,j})\}_{0\leq i,j\leq
d}\rangle.\label{rst}
\end{eqnarray}

\begin{prop}\label{ideal} Let ${}^{rs}\mathcal{T}$ be as in
{\rm(\ref{rst})}. Then ${}^{rs}\mathcal{T}$ is an ideal of
$\mathcal{T}$.
\end{prop}
\proof It suffices  to show that $L(^{rs}T_{ij})L(^{pq}T_{lm})\in
{{}^{rs}\mathcal{T}}$ and $L({}^{pq}T_{lm})L({}^{rs}T_{ij})\in
{}^{rs}\mathcal{T}.$

If $j\neq l,$ then $L(^{rs}T_{ij})L(^{pq}T_{lm})=0.$ Suppose $j=l.$
Since $H_{d-i,d-j}(d)\otimes H_{i,j}(n-d)\in M_{i,j}^{(n,d)}$ and
$H_{d-j,d-m}(d)\otimes H_{j,m}(n-d)\in M_{j,m}^{(n,d)},$ we obtain
$(H_{d-i,d-j}(d)H_{d-j,d-m}(d))\otimes (H_{i,j}(n-d)H_{j,m}(n-d))\in
M_{i,m}^{(n,d)}$. It follows that there exist   scalars
$\beta_{u,v}$ such that
\begin{eqnarray*}
&&(H_{d-i,d-j}(d)H_{d-j,d-m}(d))\otimes (H_{i,j}(n-d)H_{j,m}(n-d))\\
&=&\sum_{u=0}^{\min(i,d-i)}\sum_{v=0}^{\min(i,n-d-i)}\beta_{u,v}(E_u^{(d,d-i)}H_{d-i,d-m}(d))\otimes
(E_v^{(n-d,i)}H_{i,m}(n-d)).
\end{eqnarray*}
By Lemma \ref{change},
$$
\begin{array}{rcl}
&&^{rs}T_{ij}\; {}^{pq}T_{jm}\\
&=&\delta_{r,p}\delta_{s,q}(H_{d-i,d-j}(d)E^{(d,d-j)}_rH_{d-j,d-m}(d))\otimes
(H_{i,j}(n-d)E^{(d,j)}_sH_{j,m}(n-d)))\\
&=&\delta_{r,p}\delta_{s,q}(E^{(d,d-i)}_r\otimes
E^{(d,i)}_s)((H_{d-i,d-j}(d)H_{d-j,d-m}(d))\otimes (H_{i,j}(n-d)H_{j,m}(n-d)))\\
&=&\delta_{r,p}\delta_{s,q}\beta_{r,s}(E_r^{(d,d-i)}H_{d-i,d-m}(d))\otimes
(E_s^{(d,i)}H_{i,m}(n-d)),
\end{array}
$$
where $\delta_{r,p}$ is the Kronecker delta;  and so
\begin{eqnarray}\label{15}
{}^{rs}T_{ij}\;
{}^{pq}T_{jm}=\delta_{r,p}\delta_{s,q}\beta_{r,s}{}^{rs}T_{im}.
\end{eqnarray}
It follows that $L(^{rs}T_{ij})L(^{pq}T_{jm}) \in {}^{rs}{\mathcal
T}$. Similarly $L({}^{pq}T_{lm})L({}^{rs}T_{ij})\in
{}^{rs}\mathcal{T}.$ $\qed$

By (\ref{15}) we observe that $^{rs}{\mathcal T}^{pq}{\mathcal
T}=\{0\}$ if and only if $(r,s)\neq(p,q).$ From the construction of
$\mathcal{M}^{(n,d)},$ we have
$$
{\mathcal T}=\bigoplus_{r=0}^{\lfloor
d/2\rfloor}\bigoplus_{s=0}^{\lfloor (n-d)/2\rfloor}\ {}^{rs}
{\mathcal T}.
$$


\begin{lemma}\label{43}
Let $ {}^{rs}T_{ij}$ be as in {\rm(\ref{rsta})}. Then $
{}^{rs}T_{ij}\neq0$ if and only if
 $i,j\in\{e_{r,s},e_{r,s}+1,\ldots,e_{r,s}+d_{r,s}\}$.
\end{lemma}
\proof Note that $i,j\in\{\max(r,s),\ldots,\min(d-r,n-d-s)\}$ if and
only if $0\leq r\leq \min(i,j,d-i,d-j) $ and $0\leq s\leq
\min(i,j,n-d-i,n-d-j).$ If $r$ or $s$ does not belong to above
ranges, then $^{rs}T_{ij}=0 $ by Corollary~\ref{change}.
 Since
 \begin{equation*}\label{111}
H_{d-i,d-j}^r(d)\otimes H_{i,j}^s(n-d),\quad {r\in R(d,d-i,d-j),s\in
R(n-d,i,j)}
\end{equation*}
 is a basis of $M^{(n,d)}_{i,j},$  we have
$$
\dim(M^{(n,d)}_{i,j})=(\min(i,j,d-i,d-j)+1)\times(\min(i,j,n-d-i,n-d-j)+1).
$$
All $^{rs}T_{ij}$'s generate  $M^{(n,d)}_{i,j},$ so the desired
result follows.$\qed$

For $i, j\in \{e_{r,s},e_{r,s}+1,\ldots,e_{r,s}+d_{r,s}\}$, write

$$
n^r_{ij}=\left\{
\begin{array}{ll}
\sum\limits^{j-i}_{m=0}{d-j-m\choose d-i}p_m^{(d,d-i)}(r),& i\leq j,\\
\sum\limits^{i-j}_{m=0}{i-m\choose j}p_m^{(d,d-i)}(r),& i\geq j;
\end{array}\right.$$
$$
n^{ij}_s=\left\{
\begin{array}{ll}
\sum\limits^{j-i}_{m=0}{n-d-i-m\choose j-i-m}p_m^{(n-d,i)}(s),& i\leq j,\\
\sum\limits^{i-j}_{m=0}{i-m\choose j}p_m^{(n-d,i)}(s),&  i\geq j.
\end{array}\right.
$$

By Lemma~\ref{hh} we have $(^{rs}T_{ij})(^{rs}T_{ij})^{\rm
T}=n^r_{ij}n^{ij}_sE_r^{(d,d-i)}\otimes E_s^{(n-d,i)}\neq0.$ By
computing the trace of this matrix,
 one gets $n^r_{ij}>0$ and $n^{ij}_s>0.$
By (\ref{15}), we may assume that ${}^{rs}T_{ij}
{}^{rs}T_{jl}=\beta_{r,s}(i,j,l)\;{}^{rs}T_{il}.$ Then
$\beta_{r,s}(i,j,i)=n^r_{ij}n^{ij}_s>0.$ Taking the transpose on
both sides of above equation, we obtain
$\beta_{r,s}(i,j,l)=\beta_{r,s}(l,j,i).$ By Lemma~\ref{hh} {\rm(i)}
and Lemma~\ref{change}, we have $\beta_{r,s}(i,j,l)>0$ if $i\geq
j\geq l.$ Note that
$$
\begin{array}{ccccc} {}^{rs}T_{ij}
{}^{rs}T_{jl}{}^{rs}T_{li}=\beta_{r,s}(i,j,l)\beta_{r,s}(i,l,i)\;{}^{rs}T_{ii}
=\beta_{r,s}(j,l,i)\beta_{r,s}(i,j,i)\;{}^{rs}T_{ii}.\\
{}^{rs}T_{li}{}^{rs}T_{ij}
{}^{rs}T_{jl}=\beta_{r,s}(l,i,j)\beta_{r,s}(l,j,l)\;{}^{rs}T_{ll}
=\beta_{r,s}(i,j,l)\beta_{r,s}(l,i,l)\;{}^{rs}T_{ll}.
\end{array}
$$
Hence, we have $\beta_{r,s}(i,j,l)>0$ for any $i,j,l\in
\{\max(r,s),\ldots,\min(d-r,n-d-s)\}.$

By Lemma~\ref{hh} again,
$$
({}^{rs}T_{ij}{}^{rs}T_{jl})({}^{rs}T_{ij}{}^{rs}T_{jl})^{\rm
T}=\frac{n^r_{ij}n^{ij}_sn^r_{jl}n^{jl}_s}{n^r_{il}n^{il}_s}\;{}^{rs}T_{il}({}^{rs}T_{il})^{\rm
T}.
$$
By (\ref{15}), we have \begin{eqnarray}\label{ttt1}
{}^{rs}T_{ij}{}^{rs}T_{jl}=\sqrt\frac{n^r_{ij}n^{ij}_sn^r_{jl}n^{jl}_s}{n^r_{il}n^{il}_s}\;
{}^{rs}T_{il}.
\end{eqnarray}

Let ${\rm Mat}_{d_{r,s}+1}(\mathbb{C})$ be the algebra consisting of
 all matrices with degree $d_{r,s}+1$ whose rows and columns are
indexed by $\{e_{r,s},e_{r,s}+1,\ldots,e_{r,s}+d_{r,s}\}$. Let
$E_{i,j}$ be the matrix in ${\rm Mat}_{d_{r,s}+1}(\mathbb{C})$ whose
$(i,j)$-entry is $1$ and others are $0.$

\begin{thm}
Suppose $2d< n<3d.$ Let ${\mathcal T}$ be the Terwilliger algebra of
the Johnson scheme $J(n,d).$  Then
$$
{\mathcal T}\simeq\bigoplus_{r=0}^{\lfloor
d/2\rfloor}\bigoplus_{s=0}^{\lfloor (n-d)/2\rfloor}{\rm
Mat}_{d_{rs}+1}(\mathbb{C}).
$$
\end{thm}
\proof It suffices to prove that ${}^{rs}\mathcal{T}\simeq{\rm
Mat}_{d_{r,s}+1}(\mathbb{C}).$ Define the linear mapping $\phi$ from
${}^{rs}\mathcal{T}$ to ${\rm Mat}_{d_{r,s}+1}(\mathbb{C})$ such
that $\phi(L({}^{rs}T_{ij}))=\sqrt{n^r_{ij}n^{ij}_s} E_{i,j}.$ By
(\ref{ttt1}), we have ${}^{rs}\mathcal{T}\simeq{\rm
Mat}_{d_{r,s}+1}(\mathbb{C}).\qed$

\section{$\mathcal{T}$-algebra of $J(2d,d)$}
Let $\mathcal N $ be as in (\ref{n}). In this section we shall prove
that $\mathcal{N}$ is the Terwilliger algebra $\mathcal{T}$ of
$J(2d,d)$.

Write $H_{d-i,d-j}^r:=H_{d-i,d-j}^r(d)$ and
$H_{i,j}^s:=H_{i,j}^s(d)$ for simplicity.

\begin{thm}\label{main2} Let $\mathcal{T}$ be the Terwilliger algebra of
$J (2d, d)$ and $\mathcal{N}$ be the algebra as in {\rm(\ref{n})}.
Then $\mathcal{T} =\mathcal{N}.$
\end{thm}
\proof Since $L(N_{i,j})=L(N_{i,i}(H_{d-i,d-j}\otimes H_{i,j})),$ by
\cite[Proposition 5.2]{Levstein} it is sufficient  to prove that,
for $0\leq r\leq s\leq\min(d,d-i)$,
\begin{eqnarray}\label{prove2}
E_r^{(d,d-i)}\otimes E_s^{(d,i)}+E_s^{(d,d-i)}\otimes E_r^{(d,i)}\in
\mathcal{T}_{|\Omega_i\times \Omega_i}.
\end{eqnarray}

We shall prove (\ref{prove2}) by induction on $i$ ($i$ decreases
from $d$ to $0$). For $i=d,$ it is trivial.

\medskip
\textbf{Case 1.} $\lceil\frac{d}{2}\rceil\leq i\leq d-1$.

For $0\leq s\leq d-i-1$ and $0\leq r\leq s,$ by Lemma~\ref{llpp},
\begin{eqnarray*}
&&E_r^{(d,d-i)}\otimes E_s^{(d,i)}+E_s^{(d,d-i)}\otimes
E_r^{(d,i)}\\
&=&l^{-1}_{d,d-i-1,r}p^{-1}_{d,i+1,s}\mathcal{P}_{i+1}(E_r^{(d,d-i-1)}\otimes
E_s^{(d,i+1)}+E_s^{(d,d-i-1)}\otimes E_r^{(d,i+1)}),
\end{eqnarray*}
which implies that $E_r^{(d,d-i)}\otimes
E_s^{(d,i)}+E_s^{(d,d-i)}\otimes E_r^{(d,i)}\in
\mathcal{T}_{|\Omega_i\times\Omega_i}.$ Write
$\lambda_r=p^{(d,i)}_1(r)$. By Lemma~\ref{stru2d}, we have
$$\begin{array}{l}
\quad(A+aI)_{|\Omega_i\times\Omega_i}
-\sum\limits_{s=0}^{d-i-1}\sum\limits_{r=0}^{d-i-1}
\dfrac{\lambda_r+\lambda_s+a}{2}(E^{(d,d-i)}_r\otimes
E^{(d,i)}_s+E^{(d,d-i)}_s\otimes E^{(d,i)}_r)\\
=\sum\limits_{q=0}^{d-i-1}
(\lambda_{d-i}+\lambda_q+a)(E^{(d,d-i)}_{d-i}\otimes
E^{(d,i)}_q+E^{(d,d-i)}_q\otimes
E^{(d,i)}_{d-i})+(2\lambda_{d-i}+a)E^{(d,d-i)}_{d-i}\otimes
E^{(d,i)}_{d-i}.
\end{array}$$
Then the right hand side of the equality belongs to
$\mathcal{T}_{|\Omega_i\times\Omega_i}.$ Moreover, its coefficients
are pairwise distinct,  so (\ref{prove2}) holds.

\ \ \newline \textbf{Case 2.} $0\leq i\leq
\lceil\frac{d}{2}\rceil-1.$

By Lemma~\ref{llpp} again,  (\ref{prove2}) holds. $\qed$

Next we shall decompose $\mathcal T$ as a  direct sum of some simple
ideals.

For $0\leq r,s\leq \frac{d}{2},$ define
\begin{eqnarray*}
e_{r,s}&=&\min\{i\mid 0\neq E^{(d,d-i)}_r\otimes
E^{(d,i)}_s+E^{(d,d-i)}_s\otimes E^{(d,i)}_r\in
N_{i,i}\},\\
d_{r,s}&=&|\{i\mid 0\neq E^{(d,d-i)}_r\otimes
E^{(d,i)}_s+E^{(d,d-i)}_s\otimes E^{(d,i)}_r\in N_{i,i}\}|-1.
\end{eqnarray*}
Note that $e_{r,s}=\max(r,s)$ and $e_{r,s}+d_{r,s}=\min(d-r,d-s).$

For $r, s\in\{0,1,\ldots,\min(d-i,i)\},$   define
\begin{eqnarray}
{}^{rs}T_{ij}&=&(E^{(d,d-i)}_rH_{d-i,d-j})\otimes
(E^{(d,i)}_sH_{i,j})+(E^{(d,d-i)}_sH_{d-i,d-j})\otimes
(E^{(d,i)}_rH_{i,j}),\label{rst2a}\\
 {}^{rs}{\mathcal
T}&=&\langle\{L(^{rs}T_{ij})\}_{0\leq i,j\leq d}\rangle.\label{rst2}
\end{eqnarray}

\begin{prop}\label{ideal2}
 Let ${}^{rs}\mathcal{T}$ be as in
{\rm(\ref{rst2})}. Then ${}^{rs}\mathcal{T}$ is an ideal of
$\mathcal{T}$.
\end{prop}
\proof It suffices to show that
$L({}^{rs}{}T_{ij})L({}^{pq}{}T_{lm})\in {}^{rs}\mathcal{T}$ and
$L({}^{pq}T_{lm})L({}^{rs}T_{ij})\in {}^{rs}\mathcal{T}.$

If $j\neq l,$ then $L(^{rs}T_{ij})L(^{pq}T_{lm})=0.$ Suppose $j=l$.
Since $H_{d-i,d-j}\otimes H_{i,j}\in N_{i,j}$ and
$H_{d-j,d-m}\otimes H_{j,m}\in N_{j,m},$ we obtain
$(H_{d-i,d-j}H_{d-j,d-m})\otimes (H_{i,j}H_{j,m})\in N_{i,m}.$ It
follows that there exist scalars $\beta_{u,v}$ such that
\begin{eqnarray*}
&&(H_{d-i,d-j}H_{d-j,d-m})\otimes (H_{i,j}H_{j,m})\\
&=&\sum_{v=0}^{\min(i,d-i)}\sum_{u=0}^{v}\beta_{u,v}(E_u^{(d,d-i)}\otimes
E_v^{(d,i)}+E_v^{(d,d-i)}\otimes E_u^{(d,i)})(H_{d-i,d-m}\otimes
H_{i,m}).
\end{eqnarray*}

By Corollary~\ref{change},
\begin{eqnarray*}
&& {}^{rs}T_{ij}{}^{pq}T_{jm} \\
&=&\delta_{r,p}\delta_{s,q}(E^{(d,d-i)}_r\otimes
E^{(d,i)}_s+E^{(d,d-i)}_s\otimes
E^{(d,i)}_r)((H_{d-i,d-j}H_{d-j,d-m})\otimes (H_{i,j}H_{j,m}))\\
&=&\delta_{r,p}\delta_{s,q}\beta_{r,s}((E_r^{(d,d-i)}H_{d-i,d-m})\otimes
(E_s^{(d,i)}H_{i,m})+(E_s^{(d,d-i)}H_{d-i,d-m})\otimes
(E_r^{(d,i)}H_{i,m})),
\end{eqnarray*}
so we have
\begin{eqnarray}\label{20}
 {}^{rs}T_{ij}{}^{pq}T_{jm}=\delta_{r,p}\delta_{s,q}\beta_{r,s}{}^{rs}T_{im}.
 \end{eqnarray}
It follows that $L(^{rs}T_{ij})L(^{pq}T_{jm}) \in {}^{rs}{\mathcal
T}$. Similarly $L({}^{pq}T_{lm})L({}^{rs}T_{ij})\in
{}^{rs}\mathcal{T}.$ $\qed$

By (\ref{20}) we observe that $^{rs}{\mathcal T}^{pq}{\mathcal
T}=\{0\}$ if and only if $(r,s)\neq(p,q).$  From the construction of
$\mathcal{N},$ we have
$$
{\mathcal T}=\bigoplus_{s=0}^{\lfloor
d/2\rfloor}\bigoplus_{r=0}^{s}\ {}^{rs} {\mathcal T}.
$$


\begin{lemma}
Let $ ^{rs}T_{ij}$ be as in {\rm (\ref{rst2a})}. Then
$^{rs}T_{ij}\neq0$ if and only if
$i,j\in\{e_{r,s},\ldots,e_{r,s}+d_{r,s}\}$.
\end{lemma}
%
\proof
 The proof is similar to that of Lemma~\ref{43} and will be omitted. $\qed$

For $i, j\in \{\max(r,s),\ldots,\min(d-r,d-s)\}$, write
$$
n^{ij}_s=\left\{
\begin{array}{ll}
\sum\limits^{j-i}_{m=0}{d-i-m\choose j-i-m}p_m^{(d,i)}(s),& i\leq j,\\
\sum\limits^{i-j}_{m=0}{i-m\choose j}p_m^{(d,i)}(s),&  i\geq j.
\end{array}\right.
$$
Similar to the proof of (\ref{ttt1}), we have
\begin{eqnarray}\label{ttt2}
{}^{rs}T_{ij}{}^{rs}T_{jl}=\sqrt\frac{n^{ij}_rn^{ij}_sn^{jl}_rn^{jl}_s}{n^{il}_rn^{il}_s}
\;{}^{rs}T_{il}.
\end{eqnarray}
Let  $\phi$ be the linear mapping from ${}^{rs}\mathcal{T}$ to ${\rm
Mat}_{d_{r,s}+1}(\mathbb{C})$ satisfying
$\phi(L({}^{rs}T_{ij}))=\sqrt{n^{ij}_rn^{ij}_s} E_{i,j}.$ By
(\ref{ttt2}), we obtain the following result.

\begin{thm}
Let ${\mathcal T}$ be the Terwilliger algebra of the Johnson scheme
$J(2d,d)$. Then
$$
{\mathcal T}\simeq\bigoplus_{s=0}^{\lfloor
d/2\rfloor}\bigoplus_{r=0}^{s}{\rm Mat}_{d_{rs}+1}(\mathbb{C}).
$$
\end{thm}

\section*{Acknowledgement}
This research is supported by NSFC(11271047), Fundamental Research
Funds for the Central University of China , Secyt-UNC
2012-2013,(05-M232) and PIP-CONICET 2010-2012 (11220090100544)
Argentina.


\end{CJK*}

\end{document}